\documentclass[a4paper,10pt]{article}
\usepackage{bbm,amssymb,amsthm,amsmath,amsfonts,url,mathdots}
\usepackage[ruled,vlined,linesnumbered]{algorithm2e}
\usepackage{color,epsfig}

\def\N{\mathbb{N}}
\def\Z{\mathbb{Z}}
\def\E{\mathbb{E}}
\def\P{\mathbb{P}}

\def\1{\mathbbm{1}}
\def\A{\mathcal{A}}
\def\M{\mathcal{M}}
\def\U{\mathcal{U}}

\def\w{\underline{w}}
\def\0{\underline{0}}
\def\z{\underline{z}}
\def\hdelta{\underline{\delta}}
\def\G{G^{-\N^+}}
\def\hassuff{\succeq}
\def\depth{\mathrm{d}}
\def\F{\mathcal{F}}

\def\pcD{\overleftarrow{D}}
\DeclareMathOperator*{\racine}{\mathcal{R}}

\DeclareMathOperator*{\tree}{\mathcal{T}}

\newtheorem{lemma}{Lemma}
\newtheorem{prop}{Proposition}
\newtheorem{defi}{Definition}
\newtheorem{remark}{Remark}

\title{Perfect Simulation of Processes With Long Memory: \\ A ``Coupling Into and From The Past'' Algorithm}
\author{Aur\'elien Garivier}

\begin{document}

\maketitle

\begin{abstract}
We describe a new algorithm for the perfect simulation of variable length Markov chains and random systems with perfect connections.
This algorithm, which generalizes Propp and Wilson's simulation scheme, is based on the idea of coupling into and from the past. It improves on existing algorithms by relaxing the conditions on the kernel and by accelerating convergence, even in the simple case of finite order Markov chains.
Although chains of variable or infinite order have been widely investigated for decades, their use in applied probability, from information theory to bio-informatics and linguistics, has recently led to considerable renewed interest. 
\end{abstract}


\section{Introduction}\label{sec:intro}
Since the publication of Propp and Wilson's seminal paper~\cite{ProppWilson96ExactSampling}, perfect simulation schemes for stationary Markov chains have been developed and implemented in several fields of applied probabilities, from statistical physics to Bayesian statistics (see for example \cite{MurdochGreen98Exact} and references therein, or \cite{Haggstrom02finitechains} for an introduction).

In 2002, Comets et al.~\cite{CometsFernandezFerrari02longMem} proposed an extension to processes with long memory: they provided a perfect simulation algorithm for stationary processes called random systems with complete connections~\cite{OnicescuMihoc35achaines,OnicescuMihoc35bchaines} or \emph{chains of infinite order}~\cite{Harris55chains}. These processes are characterized by a transition kernel that specifies (given an infinite sequence of past symbols) the probability distribution of the next symbol. Following  \cite{Lalley86regen,Lalley00regen,Berbee87chains}, their work was based on the idea of exploiting the regenerative structures of these processes.
The algorithm used by Comets et al. relied on renewal properties that resulted from summable memory decay conditions. As a by-product of the simulation scheme, the authors proved the existence (and uniqueness) of a stationary process $f$ under suitable hypotheses.

However, these conditions on the kernel seemed quite restrictive and unnecessary. Gallo~\cite{Gallo10perfectSim} and Foss et al.~\cite{FossTweedie98simu} showed that different  coupling schemes could be designed under alternative assumptions that do not even require kernel continuity. 
Moreover, the coupling scheme described in~\cite{CometsFernandezFerrari02longMem} strongly relies on regeneration, not on coalescence. Contrary to Propp and Wilson's algorithm, this coupling scheme does not converge for all mixing Markov chains; when it does converge, it requires a larger number of steps. Recently, De Santis and Piccioni~\cite{DesantsiPiccioni12perfectSim} tried to combine the two algorithms by providing a hybrid method that works with two regimes: coalescence for short memory and regeneration on long scales.

This paper aims to fill the gap between long- and short-scale methods by providing a relatively elaborate, yet elegant coupling procedure that relies solely on coalescence. For a (first order) Markov chain process, this procedure is equivalent to Propp and Wilson's algorithm. However, our procedure makes it possible to manage more general, infinite memory processes characterized by a continuous transition kernel, as defined in Section~\ref{sec:notations}.

This procedure is based on the idea of exploiting \emph{local} coalescence instead of global loss of memory properties. 
From an abstract perspective, the algorithm described in Section~\ref{sec:perfectsim} simply consists in running a Markov chain on an infinite, uncountable set until the first hitting time of a given subset of states. Its concrete implementation involves a dynamical system on a set of labeled trees described in Section~\ref{sec:algo}. 

Alternatively, this algorithm may be linked to the algorithm described by Kendall~\cite{kendall98PS}, whose adaptation of Propp and Wilson's idea perfectly simulates the equilibrium distribution of a spacial birth-and-death process to precisely sample from area-interaction point processes. Analogous to the present article, Kendall's work was based on the idea that if---following the coupled transitions---all possible initial patterns at time $t<0$ lead to the same configuration at time $0$, then this configuration has the expected distribution: Whenever such a coalescence is observable, perfect simulation is possible.
This algorithm was later generalized by Wilson, who named it `coupling into and from the past' (CIAFTP, see~\cite{wilson00PSreadonce}, Section~7), an appropriate term for the algorithm described in Section~\ref{sec:algo}.

We show that this perfect simulation scheme converges under less restrictive hypotheses than previously required. The paper also provides a detailed description of finite, but large order Markov chains (or variable order Markov chains, see~\cite{buhlmann1999}) because they prove very useful in many applications (e.g. information theory \cite{rissanen1983,willems1995} or bio-informatics \cite{busch2009}): our algorithm compares favorably with Propp and Wilson's algorithm on the extended chain in terms of computational complexity; it also compares favorably with the procedure of \cite{CometsFernandezFerrari02longMem} in terms of convergence speed.

The paper is organized as follows:
Section~\ref{sec:notations} presents the notation and necessary definitions.
Section~\ref{sec:perfectsim} contains the conceptual description of the perfect simulation schemes. An update rule constructed in Section~\ref{sec:coupling} served as the key tool.
Section~\ref{sec:algo} contains the detailed description of the algorithm, and Section~\ref{sec:boundingSize} presents some  elements of the complexity analysis.  Section~\ref{sec:examples} demonstrates the relative (compared with other coupling schemes) weakness of the assumptions required for the algorithm to converge.
Finally, proofs of technical results are included in the Appendix.

\section{Notation and definitions}\label{sec:notations}
The following section introduces the notation used for algorithms, proofs, and results. Although we primarily relied on standard notation, specific symbols were required, especially concerning trees. This somewhat unusual notation, which is necessary to expose the algorithm as clearly as possible, is central to this paper.
\subsection{Histories}
Corresponding to \cite{CometsFernandezFerrari02longMem}, $G$ denotes a finite alphabet whose size is denoted by $|G|$.
For $k\in\N$, $G^{-k}$ denotes the set of all sequences $(w_{-k}, \dots, w_{-1})$, and $G^* = \cup_{k\geq 0} G^{-k}$.
By convention, $\varepsilon$ denotes the empty sequence, and $G^0 = \{\varepsilon\}$.
Referred to as the space of \emph{histories}, the set of $G$-valued sequences indexed by the set of negative integers is denoted by $\G$. 
For $-\infty\leq a\leq b <0$ and $w\in \G$, the sequence $(w_a,\dots,w_b)$ is denoted by $w_{a:b}$. An element $w_{-\infty:-1}\in\G$ is denoted by $\w$. 
For $w\in G^{-k}$, we write $|w|=k$; for $\w\in \G$,  we write $|\w|=\infty$.
For every negative integer $n$, we define the projection $\Pi^{n}: \G \to G^{n}$ by $\Pi^n(\w) = w_{n:-1}$.

\begin{figure}[htbp]
\begin{center}
\includegraphics[width=10cm]{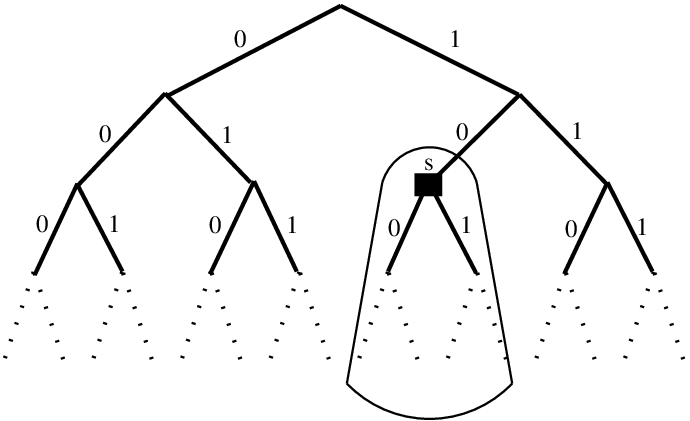}
\caption{Trie representation of $\G$ for the binary alphabet $G=\{0,1\}$. The square represents $s=(0,1)\in G^{-2}$, and $\tree(s)$ is circled.
Note that the symbols $(s_{-2}, s_{-1})=(0,1)$ are to be read \emph{from bottom} (node $s$) \emph{to top} (root).}
\label{fig:completeTree}
\end{center}
\end{figure}

A  \emph{trie} is a rooted tree whose edges are labeled with elements of $G$. 
An element $\w\in\G$ can be represented by a path in the infinite, complete trie starting from the root and successively following the edges labeled as  $w_{-1}, w_{-2}, \dots$
A finite sequence $s\in G^*$ is represented by an internal node of this infinite trie. 
Figure~\ref{fig:completeTree} illustrates this representation for the binary alphabet $G=\{0,1\}$.

\subsection{Concatenation and suffix}
For the two sequences $w_{a:b}$ and $z_{c:d}$, with $-\infty\leq a\leq b<0$ and $-\infty<c\leq d<0$, the concatenation of $w_{a:b}$ and $z_{c:d}$ is written $w_{a:b}z_{c:d} = (w_a,\dots,w_b,z_c,\dots,z_d)$. In particular, if we take $a=-\infty$, then this defines the concatenation $\z s$ of a history $\z$ and an $n$-tuple $s\in G^{|s|}$. Note that this notation is different from the convention taken in  \cite{CometsFernandezFerrari02longMem}.
If $a>b$, then $w_{a:b}$ is the empty sequence  $\varepsilon$.

Let $h\in G^*\cup \G$. If $s\in G^*$ is such that $|h|\geq |s|$ and $h_{-|s|:-1} = s$, we say that $s$ is a \emph{suffix} of $h$ and we write $h \hassuff s$. This defines a partial order $\hassuff$ on $G^*\cup \G$.

\subsection{Metric}
Equipped with the product topology and with the ultra-metric distance $\delta$ defined by \[\delta(\w, \z) = 2^{\sup\{k<0 : w_k\neq z_k\}},\] $\G$ is a complete and compact set.
A ball $B \subset\G$ is a set $\left\{\z s : \z\in\G\right\}$ for some $s\in G^*$. In reference to the trie representation of $\G$, we write $s=\racine(B)$ for the \emph{root} of $B$, and $\tree(s)=B$ for the \emph{tail} of $s$ (see Figure~\ref{fig:completeTree}). 
Note that $\tree(\varepsilon) = \G$.

The set of probability distributions on $G$ is denoted by $\M(G)$; it is endowed with the total variation distance
\[|p-q|_{TV} = \frac{1}{2}\sum_{a\in G} |p(a)-q(a)| = 1-\sum_{a\in G} p(a)\wedge q(a)\;,\]
where $x\wedge y$ is the minimum of $x$ and $y$.

\subsection{Complete suffix dictionaries}
\begin{figure}[htbp]
\begin{center}
\input{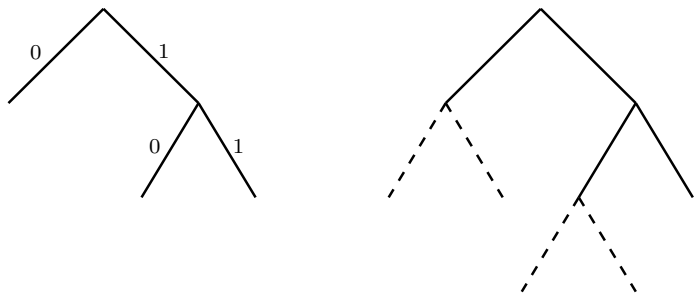}
\caption{Trie representation of a CSD for the binary alphabet $G=\{0,1\}$. Left: the trie representing the complete suffix dictionary $D = \{0, 01, 11\}$. Right: $\{00, 10, 001, 101, 11\}\hassuff \{0, 01, 11\}$.}
\label{fig:exampleCSD}
\end{center}
\end{figure}

A (finite or infinite) set $D$ of elements of $G^*$ is called a \emph{complete suffix dictionary} (CSD) if one of the following equivalent properties is satisfied:
\begin{itemize}
 \item Every sequence $\w\in\G$ has a unique suffix in $D$: \[\forall\ \w\in\G, \, \exists ! s\in D : \w \hassuff s\;.\]
 \item $\big\{\tree(s) : s\in D\big\}$ is a partition of $\G$, in which case we write \[ \G = \bigsqcup_{s\in D} \tree(s)\;.\]
\end{itemize}
A CSD can be represented by a trie, as illustrated in Fig.~\ref{fig:exampleCSD}.
This representation suggests that the \emph{depth} of CSD $D$ is defined as the depth of this trie: 
\[\depth(D) = \sup\big\{|s| : s\in D\big\}\;.\]
Note that $\depth(D)=+\infty$  if $D$ is infinite. 
The smallest possible CSD is $\{\epsilon\}$ (its trie is reduced to the root): it has a depth of $0$ and a size of $1$. The second smallest is $G$ with a depth of $1$.
 
If a finite word $h\in G^*$ has a (unique) suffix in $D$, we write $h\hassuff D$.
If $D$ and $D'$ are two CSDs such that $s\hassuff D$ as soon as $s\in D'$, we write $D'\hassuff D$. This means that the trie representing $D'$ entirely covers that of $D$, as illustrated in Fig.~\ref{fig:exampleCSD}.

\subsection{Piecewise constant mappings}\label{subsec:PCmappings}
For a given CSD $D$, we say that a mapping $f$ defined on $\G$ is \emph{$D$-constant} if,
\[\forall s\in D ,\, \forall \w, \z\in \tree(s), f(\w)=f(\z) \;.\]
The mapping $f$ is constant if and only if it is $\{\epsilon\}$-constant, and $f$ is called \emph{piecewise constant} if there exists a CSD $D$ such that $f$ is $D$-constant.
For every $h\in G^*$ we define 
\[f(h) = f\big(\tree(h)\big) = \left\{ f(z) : z\in \tree(h)\right\}\;.\]
Note that, by definition, $f(h)$ is a set. 
However, if $f$ is $D$-constant and if $h\hassuff D$, then $f(h)$ is a singleton (a set containing exactly one element).

\begin{figure}[htbp]
\begin{center}
\input{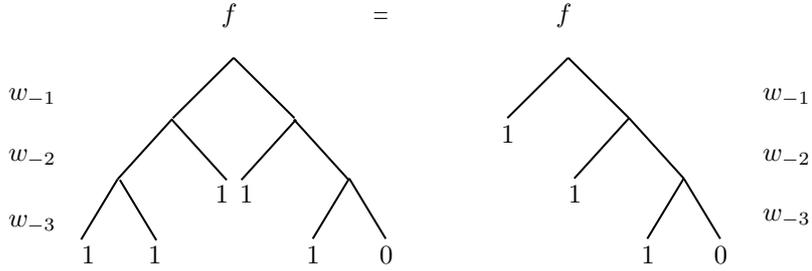}
\caption{Two representations as labeled tries of the piecewise constant function $f$ defined for the binary alphabet $\{0,1\}^{-\N^+}$ by $f(\w) = 0$ if $\w \in \tree(111)$, and $f(\w) = 1$ otherwise. In the second representation, the trie is minimal: it has been obtained from the first trie by recursively pruning leaves with identical images.}
\label{fig:exPCfunction}
\end{center}
\end{figure}

Let $f$ be a piecewise constant mapping. The set of all CSDs such that $f$ is $D$-constant has a minimal element when ordered by the inclusion relation: $D^f$ denotes the \emph{minimal CSD} of $f$.
The minimal CSD $D^f$ is such that, if $s\in D^f$, there exists $w\in G^*$ such that $s'=ws_{-|s|+1:-1}\in D^f$ and $f(s)\neq f(s')$.
If $f$ is $D$-constant, then $D^f$ can be obtained by recursive pruning of $D$, that is by pruning the nodes whose children are leaves with the same $f$ value (repeating this operation for as long as possible).
A $D$-constant mapping $f$ can be represented by the trie $D$ if each leaf $s$ of $D$ is labeled with the common value of $f(\w)$ for $w\in\tree(s)$. 
Figure~\ref{fig:exPCfunction} illustrates the trie representation of a piecewise constant function as well as the pruning operation.

\subsection{Probability transition kernels}
A mapping $P: \G \to \M(G)$ is called a \emph{probability transition kernel}, and we write $P(\cdot|\w)$ for the image of  $\w \in \G$.
We say that $P$ is \emph{continuous} if it is continuous as an application from $\left(\G, \delta\right)$ to $\left(\M(G), |\cdot|_{TV}\right)$.
For $s \in G^*$, we define the \emph{oscillation} of $P$ on the ball $\tree(s)$ as
\[\eta_P(s) = \sup\Big\{\,\big|P(\cdot|\w) - P(\cdot|\z) \big|_{TV}  : \w, \z \in \tree(s) \,\Big\}.\]

We say that a process $(X_t)_{t\in \Z}$  with a distribution $\nu$ on $G^{\Z}$ (equipped with the product topology and the product sigma-algebra) is \emph{compatible} with kernel $P$ if the latter is a version of the one-sided conditional probabilities of the former; that is 
\[\nu\left(X_i = g | X_{i+j}=w_j \;\text{for all}\ j\in-\N^+ \right) = P(g|\w)\]
for all $i\in\Z, g\in G$, and $\nu$-almost every $\w$.
A classical but key remark states that $S_t = (\dots, X_{t-1}, X_t) , t\in\Z$,
is  a homogeneous Markov chain on the compact ultra-metric state space $\G$ with a transition kernel $Q$ given by the relation,
\[\forall \w , \z\in \G, \quad Q(\z | \w) = P(z_{-1} | \w) \1_{\bigcap_{i<0}\{z_{i-1} = w_i\}}\;.\]

\subsection{Update rules}\label{sec:couplings}
An application $\phi: [0,1[\times \G\to G$ is called an \emph{update rule} for a kernel $P$ if, for all $\underline{w}\in G^{-\N^+}$ and for all $g\in G$, the Lebesgue measure of $\{u\in[0,1[: \phi(u, \underline{w}) = g\}$ is equal to $P(g|\underline{w})$. 
In other words, if $U$ is a random variable uniformly distributed on $[0,1[$, then $\phi(U, \w)$ has a distribution $P(\cdot|\w)$ for all $\w\in \G$.
For any continuous kernel $P$, Section~\ref{sec:coupling} details the construction of an update rule $\phi_P$ such that,
\begin{equation}\label{eq:coupling}
\text{for all}\ s\in G^*, \;0\leq u < 1-|G|\eta_P(s) \implies \phi_P(u,\cdot) \text{ is constant on }\tree(s)\;.
\end{equation}
The following lemma (proved in the Appendix) states the basic observation that permits us to design an algorithm working in finite time, which is applicable even for kernels that are not piecewise continuous. 
\begin{lemma}\label{lem:phiPC}
For all $u\in[0,1[$, the mapping $\w\to\phi_P(u, \w)$ is continuous, i.e. piecewise constant.
\end{lemma}

\section{Abstract description of the perfect simulation scheme}\label{sec:perfectsim}
Given a continuous transition kernel $P$, two questions arise:
\begin{enumerate}
 \item Does a stationary distribution $\nu$ compatible with $P$ exist? If it exists, is it unique?
 \item If $\nu$ exists, how can we sample finite trajectories from that distribution?
\end{enumerate}
Several authors \cite{CometsFernandezFerrari02longMem,DesantsiPiccioni12perfectSim,Gallo10perfectSim} have contributed to answering these questions in the past decade. Their approach consisted in showing that there exists a simulation scheme that draws samples of $\nu$. This algorithm was based on the idea of coupling from the past. In accordance with these authors, we addressed these questions by constructing a new perfect simulation scheme that requires looser conditions on the kernel and that converges faster than existing algorithms. The following section describes the general principle of this algorithm. Practical details concerning its implementation are given in Section~\ref{sec:algo}.

\subsection{Perfect simulation by coupling into the past}

\begin{figure}[htbp]
\begin{center}
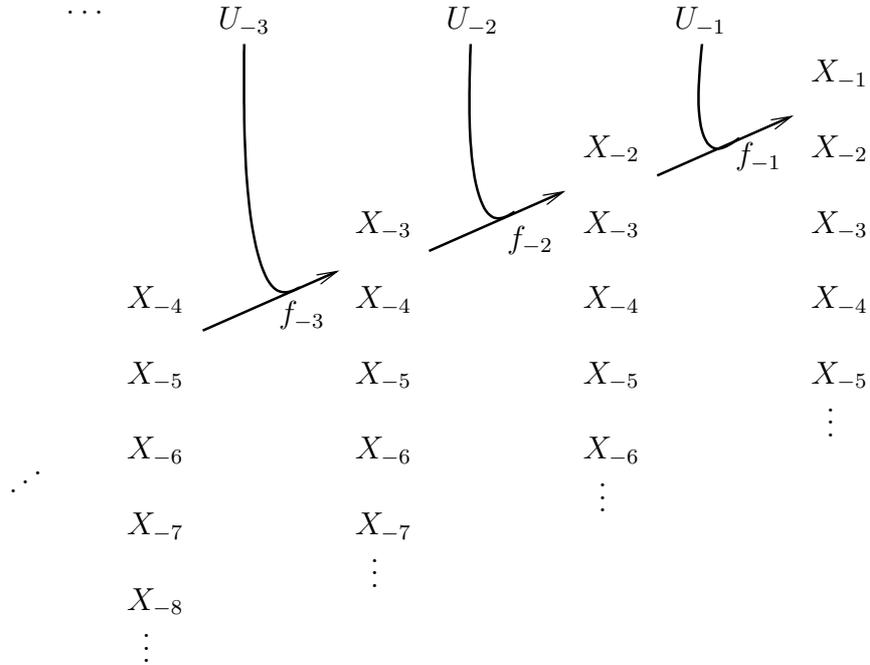
\caption{Perfect simulation scheme.}
\label{fig:principle}
\end{center}
\end{figure}

Let $n$ be a negative integer. In order to draw $(X_{n}, \dots, X_{-1})$ from a stationary distribution compatible with $P$, we use a semi-infinite sequence of independent random variables $(U_t)_{t<0}$ defined on a probability space $(\Omega, \A, P)$ and uniformly distributed on $[0,1[$. The variable $X_t$ is deduced from $U_t$ and from the past symbols $X_{t-1},X_{t-2},\dots$, as depicted in Fig.~\ref{fig:principle}. Those past symbols are unknown, but the regularity of $P$ makes it nevertheless possible to sometimes compute $X_t$.

For each $t<0$, let $f_t$ be the random function $\G\to \G$ defined by $f_t(\w) = \w\phi_P(U_t, \w)$.\footnote{Regarding measurability issues: if the set of functions $\G\to \G$ is equipped with the topology induced by the distance $\hdelta$ defined by 
\[\hdelta(f_1,f_2) = \sum_{w \in G^*} (2|G|)^{-|w|} \sup\Big\{\delta\big(f_1(\z_1), f_2(\z_2)\big): \z_1,\z_2\in\tree(w)\Big\} \;\]
and with the corresponding Borel sigma-algebra, then the measurability of $f_t$ follows from Lemma~\ref{lem:phiPC}.
} Beware of the index shift: if $\z=f_t(\w)$, then $z_{-1} = \phi_P(U_t, \w)$ and $z_i = w_{i+1}$ for $i<-1$.

In addition, let $F_t = f_{-1}\circ\dots\circ f_t$ and, for any negative integer $n$, $H^n_{t} = \Pi^n \circ F_{t}$.
Proposition~\ref{lem:phiPC} below shows that the continuity of $P$ implies that $H^n_{t}$ is piecewise constant. 
We define \[\tau(n) = \sup\{t <0 : H^n_{t} \hbox{ is constant}\} \;, \]
where, by convention, $\tau(n)=-\infty$ if $H^n_{t}$ is not constant for all $t<-1$.
When $\tau(n)$ is finite, the result of the procedure is the image $\{X_{n:-1}\}$ of the constant mapping $H_{\tau(n)}^n$. We can easily verify that $X_{n:-1}$ has the expected distribution (see~\cite{ProppWilson96ExactSampling,CometsFernandezFerrari02longMem}).
\begin{remark}\label{rem:cst}
For $t>n$, $H_t^n$ cannot be constant because it holds that $\left(H_t^n(\w) \right)_k = w_{k}$ for all $n\leq k <t$.
Thus, $\tau(n) =  \sup\{t \leq n : H^n_{t} \hbox{ is constant}\} \leq n$.
\end{remark}
Observe also that the sequence $(\tau(n))_n$ is a non-increasing sequence of stopping times with respect to the filtration $(\F_s)_s$, where $\F_s = \sigma(U_t: t\geq s)$ when $s$ decreases.

From a theoretical perspective, this CIAFTP algorithm simply consists in running an instrumental Markov chain until a given hitting time. 
In fact, the recursive definition given above shows that the sequence $(H_t^n)_{t\leq 1}$ is a homogeneous Markov chain on the set of functions $\G\to G^n$.
The algorithm terminates when this Markov chain hits the set of constant mappings.
Such a procedure seems to be purely abstract because it involves infinite, uncountable objects. 
However, Section~\ref{sec:algo} shows how this Markov chain on the set of functions $\G\to G^n$ can be handled with a finite memory. Before we provide the detailed implementation of the algorithm, we first present the construction of the update rule and the sufficient conditions for the finiteness of the stopping time $\tau(n)$ in Section~\ref{sec:coupling}.

\subsection{Constructing the update rule $\phi_P$}\label{sec:coupling}
The algorithm that is abstractly depicted above and detailed in Section~\ref{sec:algo} crucially relies on the update rule $\phi_P$ that satisfies Equation~\eqref{eq:coupling}. We present here the construction of this update rule for a given continuous kernel $P$. 
In short, for each $k$-tuple $z\in G^{-k}$, the construction of $\phi_P$ relies on a coupling of the conditional distributions $\big\{P(\cdot | \z) : \z \in \tree(z)\big\}$.
The simultaneous construction of all of these couplings requires a few definitions and properties, which we state in this section and prove in the Appendix. 

Provide $G$ with any order $<$, so that $\G$ can be equipped with the corresponding lexicographic order : $\w<\z$ if there exists $k\in-\N$ such that, for all $j>k, w_j=z_j$, and $w_k<z_k$.
The continuity of $P$ is locally quantified by some coupling factors, which we define along with the coefficients that are necessary for the construction of the update rule $\phi_P$. For all $g\in G$, let $A_{-1}(\varepsilon)=a_{-1}(g | \varepsilon) = 0$. For all $k\in\N$ and all $z\in G^{-k}$, let 
\begin{align}
a_k(g | z_{-k:-1}) &=  \inf\big\{P(g|\w): \w\in \tree\left( z_{-k:-1} \right)\big\} \;,\nonumber\\
A_k(z_{-k:-1}) &= \sum_{g\in G}a_k(g|z_{-k:-1}) \;,\nonumber\\
A^-_k& = \inf_{s\in G^{-k}} A_k(s) \;,\nonumber\\
\alpha_{k}(g|z_{-k:-1}) &= A_{k-1}(z_{-k+1:-1}) + \sum_{h<g} \left\{a_{k}(h | z_{-k:-1}) - a_{k-1}(h | z_{-k+1:-1})\right\} \;,\label{eq:defalpha}\\
\beta_{k}(g|z_{-k:-1}) &=  A_{k-1}(z_{-k+1:-1}) + \sum_{h\leq g} \left\{a_{k}(h | z_{-k:-1}) - a_{k-1}(h | z_{-k+1:-1})\right\} \;.\label{eq:defbeta}
\end{align} 
Note that, with our conventions, $a_0(g|\varepsilon) = \inf\{P(g|\z) : \z\in\G\}$. Moreover, if $s$ and $s'\in G^*$ are such that $s\hassuff s'$, then for all $g\in G$ it holds that $a_k(g | s)\geq a_k(g | s')$, $A_k(s) \geq A_k(s')$, and the sequence $\left(A^-_k\right)_k$ is non-decreasing.

The following propositions gather some elementary ideas that are used in the following; they are proved in the Appendix.
\begin{prop}\label{prop:etavsA}
The coupling factors of the kernel $P$ satisfy the following inequalities: for all $s\in G^*$,
 \begin{equation}\label{eq:etavsA}
 1-(|G|-1)\eta_P(s) \leq A_{|s|}(s) \leq 1-\eta_P(s)\; .
\end{equation}
\end{prop}
\begin{prop}\label{prop:Pcont}
The following assertions are equivalent:
\begin{enumerate}
 \item[(i)] the kernel $P$ is continuous;
 \item[(ii)] for every $\w\in\G$, $\eta_P(w_{-k:-1})$ tends to $0$ when $k$ goes to infinity;  
 \item[(iii)] when $k$ goes to infinity, \[\sup\big\{\eta_P(s) : s\in G^{-k}\big\}\to 0\;;\]
 \item[(iv)] for all $\w\in\G$, $A_k(w_{-k:-1}) \to 1$ when $k\to \infty$;
 \item[(v)] $A^-_k\to 1$ when $k$ goes to infinity.
\end{enumerate}
\end{prop}
\begin{prop}\label{prop:partition}
Let $P$ be a continuous kernel, and let $\alpha_k(\cdot|\cdot)$ and $\beta_k(\cdot|\cdot)$ be defined as in~\eqref{eq:defalpha} and~\eqref{eq:defbeta}. Then, for every $\w\in \G$,
\[[0,1[ = \bigsqcup_{g\in G, k\in\N}[\alpha_k(g|w_{-k:-1}), \beta_k(g|w_{-k:-1}) [\;.\]
In other words, for every $u\in[0,1[$ and every $\w\in\G$, there exists a unique $k\in\N$ and a unique $g\in G$ such that $u\in\left[\alpha_k(g|w_{-k:-1}), \beta_k(g|w_{-k:-1})\right[$.
\end{prop}

\begin{figure}[htbp]
\begin{center}
\includegraphics[width=12cm]{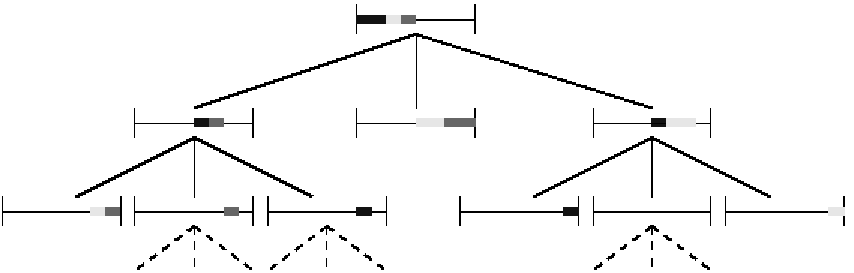}
\caption{
Graphical representation of an update rule $\phi_P$ on the alphabet $\{0,1,2\}$: for each $w_{-k:-1}$, the intervals $[\alpha_k(g|w_{-k:-1}), \beta_k(g|w_{-k:-1})[$ are represented in black ($g=0$), light grey ($g=1$), and medium grey ($g=2$).
For example, $P(1|1) = \alpha_0(1|\varepsilon) + \alpha_1(1|1) = 1/8+1/4$, and $P(0|00) = \alpha_0(0|\varepsilon) + \alpha_1(0|0) + \alpha_2(0|00) = 1/4+1/8+0$. 
}
\label{fig:coupling}
\end{center}
\end{figure}

Figure~\ref{fig:coupling} illustrates Proposition~\ref{prop:partition} on a three-symbol alphabet.
Thanks to Proposition~\ref{prop:partition}, we can now define the following update rule and verify that it satisfies  Equation~\eqref{eq:coupling}.
\begin{defi}\label{def:phi}
Let  $\phi_P:[0,1[\times\G\to G$ be defined as 
\[\phi_P(u, \w) = \sum_{g\in G, k\in\N} g\1_{[\alpha_k(g|w_{-k:-1}), \beta_k(g|w_{-k:-1})[}(u)\;.\]
In words, for every $u\in[0,1[$ and for every $\w\in\G$, $\phi_P(u, \w)$ is the unique symbol $g\in G$ such that there exists $k\in\N$ satisfying
\[u\in\big[\alpha_k(g|w_{-k:-1}), \beta_k(g|w_{-k:-1})\big[\;.\]
\end{defi}
\begin{prop}\label{prop:phiOk}
The mapping $\phi_P$ of Definition~\ref{def:phi} is an update rule satisfying Equation~\eqref{eq:coupling}:
\[\forall s\in G^*,\, \forall u\in[0,1],\, \forall \w, \z \in \tree(s),\; \quad u < A_{|s|}(s) \implies \phi_P(u,\w)=\phi_P(u,\z)\;.\]
As a consequence, for every $s\in G^*$ and every $u< A_{|s|}(s)$, we can define $\phi_P(u, s)$ as the value common to $\phi_P(u, \w)$ for all $\w\in \tree(s)$.
\end{prop}

\subsection{Convergence}
Sufficient conditions on $P$ are given in \cite{CometsFernandezFerrari02longMem,DesantsiPiccioni12perfectSim} to ensure that $\tau(n)$ is almost surely finite (or even that $\tau(n)$ has bounded expectation).
In addition, the authors prove that the almost-sure finiteness of $\tau(n)$ is a sufficient condition to prove the existence and uniqueness of a stationary distribution $\nu$ compatible with $P$ (see \cite{CometsFernandezFerrari02longMem}, Theorem 4.3 and corollaries 4.12 and 4.13). As a by-product, they obtained a simulation algorithm for sample paths of $\nu$: if $U = (U_t)_{t\in\Z}$ is a sequence of independent, uniformly distributed random variables, then $\Phi:[0,1]^\Z\to G^\Z$ 
 can be defined such that  $\Phi(U)_t = \phi_P(U_{t-1}, \Phi(U)_{-\infty:t-1})$ for all $t$, and 
\[\nu = \P(\Phi(U)\in\cdot),\]
the law of $\Phi(U)$, is stationary and compatible with $P$.

However, in \cite{CometsFernandezFerrari02longMem}, the authors impose fairly restrictive conditions on $P$: they require that
\[\sum_{m\geq 0}^\infty\prod_{k=0}^m A_k^- = \infty\;,\]
and, in particular, that the chain satisfies the Harris condition \[A_0^- =\sum_{g\in G}\inf \left\{P(g|\z): \z\in\G\right\}>0\;.\] By using a Kalikow-type decomposition of the kernel $P$ as a mixture of Markov chains of all orders, the authors prove that the process regenerates and that the stopping time 
\[\tau'(n)= \sup\{t \leq n : H^t_{t} \hbox{ is constant}\}\]
is almost surely finite under these conditions.

This condition is clearly sufficient but certainly not necessary for $\tau(n)$ to be finite. Consider, for example, a first-order Markov chain: Although Propp and Wilson \cite{ProppWilson96ExactSampling} have shown that the stopping time $\tau(n)$ of the optimal update rule is almost surely finite for every mixing chain (and, under some conditions, that $\tau(n)$ has the same order of magnitude as the mixing time of the chain), $\tau'(n)$ is almost surely infinite as soon as the Dobrushin coefficient $A_0^-$ of the chain is $0$.
In this paper, we close the gap by providing a Propp--Wilson procedure for general continuous kernels that may converge even if the process is not regenerating. For first-order Markov chains,  $\phi_P(u, \w)$ depends only on $w_{-1}$, and the algorithm presented in this paper corresponds to Propp and Wilson's exact sampling procedure.

Since the publication of \cite{CometsFernandezFerrari02longMem}, these results have been generalized    \cite{Gallo10perfectSim,DesantsiPiccioni12perfectSim}, which included relaxing the conditions on the kernel and proposing other particular conditions for different cases. Gallo \cite{Gallo10perfectSim} showed that the kernel $P$ need not be continuous to ensure the existence of $\nu$, nor to ensure the finiteness of $\tau(n)$: he gives an example of a non-continuous regenerating chain (see also the final remark of Section~\ref{sec:examples}).
De Santis and Piccioni  \cite{DesantsiPiccioni12perfectSim} have proposed another algorithm, which combines the ideas of \cite{CometsFernandezFerrari02longMem} and \cite{ProppWilson96ExactSampling}: they  propose a hybrid simulation scheme that works with a Markov regime and a long-memory regime.
We take a different, more general approach. Our procedure generalizes the sampling schemes of \cite{CometsFernandezFerrari02longMem} and \cite{ProppWilson96ExactSampling} in a single, unified framework.

\section{The coupling into and from the past algorithm}\label{sec:algo}
This section gives a detailed description of the algorithm that permits us to effectively compute the mappings $H_t^n$.  
The problem is that the mapping domain is the infinite space $\G$, so that no naive implementation is possible. We are able to solve this problem because, for each $t$, the mapping $H_t^n$ is piecewise continuous and thus can be represented by a random but finite object: namely, by its trie representation defined in Section~\ref{subsec:PCmappings}.

\subsection{Description of the algorithm}
Consider a continuous kernel $P$ and its update rule $\phi_P$ given by Definition~\ref{def:phi}.
For each $u\in[0,1[$, Proposition~\ref{lem:phiPC} shows that the mapping $\phi_P(u, \cdot)$ is piecewise constant; its minimal CSD is written $D(u) = D^{\phi_P(u,\cdot)}$. 
Algorithm~\ref{algo:rec} shows how the mappings $H_t^n$ (defined in Section~\ref{sec:perfectsim}) can be constructed recursively using only finite memory. For simplicity, this algorithm is presented as a pseudo-code involving mathematical operations and ignoring specific data structures. It should, however, be easy to deduce a real implementation from this pseudo-code. A matlab implementation is available online at \url{http://www.math.univ-toulouse.fr/~agarivie/Telecom/context/}. It contains a demonstration script illustrating the perfect simulation of the processes mentioned in Sections \ref{subsec:finiteCT} and \ref{sec:examples}.

\begin{algorithm}
 \caption{Coupling from and into the past for continuous kernels. \label{algo:rec}}
  \BlankLine
  \KwIn{update rule $\phi_P$, size $-n$ of the path to sample}
  \BlankLine
  $t \gets 0$\;
  $D_t^n \gets G^{n}$\;
  for all $s\in G^n, H_t^n(s) \gets \{s\}$\;
  \While{$|D_{t}^n|>1$}{
    $t\gets t-1$\;
    \textbf{draw} $U_t \sim\U([0,1[)$ independently\;
    $D(U_t) \gets$ the minimal trie of $U_t$\;
    \ForEach{$s\in D(U_t)$}{
      $\big\{g_t[s]\big\} \gets \phi_P(U_t,s)$\;
      \eIf{$sg_t[s]\hassuff D_{t+1}^n$}{
        $E_t^n[s] \gets \{s\}$\;
      }
      {
	$E_t^n[s] \gets \left\{h\in G^* : hg_t[s] \in D_{t+1}^n\big(sg_t[s]\big) \right\}$\;
      }
    }
    $E_t^n \gets \displaystyle{\bigcup_{s \in D(U_t)}E_t^n[s]}$\;
    \textbf{Claim 1: } $E_t^n$ is a CSD\;
    \textbf{Claim 2: } $H_t^n$ is $E_t^n$-constant, and, for all $s\in E_t^n$, $H_t^n(s) = H_{t+1}^n\big(sg_t[s]\big)$ is a singleton\;
    $D_t^n \gets$ the minimal CSD of $H_t^n$ obtained by pruning $E_t^n$
  }
  \BlankLine
  \KwOut{ $X_{n:-1}$ such that, for all $\z\in\G, \;H_t^n(\z) = \{X_{n:-1}\}$}
\end{algorithm}

For every $t<0$, the mapping  $H_t^n$ being  piecewise constant, we can define $D_t^n = D^{H_t^n}$.
Note that the definition of  $H_0^n$ in the initialization step is consistent with the general definition $H_t^n = \Pi^n\circ F_t$ because the natural definition of $F_0$ is the identity map on $\G$.
Algorithm~\ref{algo:rec} successively computes $H_{-1}^{n}, H_{-2}^{n}, \dots$ and stops for the first $t\leq n$ such that $H_t^n$ is constant.

\begin{figure}[htbp]
\begin{center}
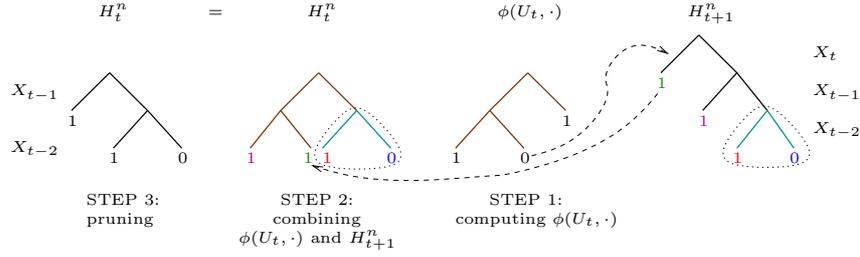
\caption{
Three-step process for the trie representation of $H_t^n$ deduced from the representations of of $H_{t+1}^n$ and $U_t$. 
Here, $G=\{0,1\}$ and $n=1$. Recall that $\phi(U_t, \cdot)$ gives $X_t$ according to $(X_{t-1}, X_{t-2},\dots)$, and that $H_{t+1}^n$ gives $X_{-1}$ according to $(X_{t}, X_{t-1}, X_{t-2}\dots)$. In Step 2, dashed lines illustrate the first case (if $(X_{t-2}, X_{t-1}) = (0,1)$, then $X_t= \phi(U_t, \dots01) = 0$ and $X_1 = H_{t+1}^n(\dots010) = H_{t+1}^n(\dots0) = 1$); for the second case, a dotted line circles the subtree of $D_{t+1}^n$ to be inserted into $E_t^n$.
}
\label{fig:exampleIter}
\end{center}
\end{figure}

The key step is the derivation of $H_t^n$ and $D_t^n$ from $H_{t+1}^n$, $D_{t+1}^n$, and $U_t$. This derivation is illustrated in Figure~\ref{fig:exampleIter} and consists of three steps:
\begin{description}
\item[STEP 1:] Compute the minimal trie $D(U_t)$ of $\phi(U_t, \cdot)$.
\item[STEP 2:] Compute the trie $E_{t}^n$ such that $H_t^n$ is $E_{t}^n$-constant by completing $D(U_t)$ with portions of $D_{t+1}^n$. Namely, for every $s\in D(U_t)$, there are two cases:
\begin{itemize}
 \item Either $sg_t[s]\hassuff D_{t+1}^n$, then knowing that $(X_{t-|s|}, \dots,X_{t-1}) = s$ together with $U_t$ and $H_{t}^n$ is sufficient to determine $X_{n:-1}$ (see the dashed lines in Figure~\ref{fig:exampleIter}),
 \item or some additional symbols in the past are required by $H_{t+1}^n$, and a subtree of $D_{t+1}^n$ has to be inserted instead of $s$ (see the dotted circled subtree in Figure~\ref{fig:exampleIter}).
\end{itemize}
\item[STEP 3:] Prune $E_{t}^n$  to obtain the minimal trie $D_{t}^n$ of $H_t^n$.
\end{description}

From a mathematical perspective, Algorithm~\ref{algo:rec} can be considered a run of an instrumental, homogeneous Markov chain on the set of $|G|$-ary trees whose leaves are labeled as $G^n$, which is stopped as soon as a tree of depth $0$ is reached. Figure~\ref{fig:exampleIter} illustrates one iteration of this chain, corresponding to one loop of the algorithm.

Algorithm~\ref{algo:rec} thus closely resembles the high-level method termed `coupling into and from the past' in~\cite{wilson00PSreadonce} (see Section~7, in particular Figure~7). In fact, in addition to coupling the trajectories starting \emph{from} all possible states at past time $t$, we use a coupling of the conditional distributions before time $t$ (that is, into the past). Our algorithm is slightly different because we want to sample $X_{n:-1}$ in addition to $X_{-1}$. The CSD $D_t^n$ corresponds to the state denoted by $X$ in~\cite{wilson00PSreadonce}, and $H_t^n$ corresponds to $F$.

\subsection{Correctness of Algorithm~\ref{algo:rec}}
To prove the correctness of Algorithm~\ref{algo:rec}, that is, the correctness of the update rule deriving $H^n_t$ from  $H_{t+1}^n$, we must verify the two claims on lines 15 and 16.

\bigskip\noindent\textbf{Claim 1: $E_t^n$ is a CSD.}\\
 Every $h\in E_t^n[s]$ is such that $h\hassuff s$. Let  $\w\in \tree(s)$, and let $\big\{h g_t[s]\big\} = D_{t+1}^n\big(\w g_t[s]\big)$. Then $h\in E_t^n[s]$ and $\w\in \tree(h)$, so that 
 \[\tree(s) = \bigsqcup_{h \in E_t^n[s]} \tree(h)\;.\]
Because $\G = \sqcup_{s\in D(U_t)} \tree(s)$, the result follows.

\bigskip
\noindent\textbf{Claim 2: $H_t^n$ is $E_t^n$-constant, and $H_t^n(s) = H_{t+1}^n\big(sg_t[s]\big)$ is a singleton for all $s\in E_t^n$.}\\
We prove that $H_t^n$ is $E_t^n$-constant by induction on $t$, and the formula for $H_t^n(s)$ is a by-product of the proof. 
For $t=0$, this is obvious if we write $E_0^n = D_0^n = G^n$.
For $t<0$, let $h\in E_t^n$. By construction, $h \hassuff D(U_t)$: denote by $s$ the suffix of $h$ in $D(U_t)$.
Then $\phi_P(U_t, h)$ is the singleton $\{g_t[s]\}$. By construction, $hg_t[s]\hassuff D_{t+1}^n$, and therefore $H_t^n(h) = H_{t+1}^n(hg_t[s])$ is a singleton by the induction hypothesis.

\subsection{Computational complexity}
For a given kernel, the random number of elementary operations performed by a computer during a run of Algorithm~\ref{algo:rec} is a complicated variable to analyze because it depends not only on the number $\tau(n)$ of iterations, but also on the size of the trees $D_t^n$ involved. Moreover, the basic operations on trees (traversal, lookup, node insertion or deletion, etc.) performed at each step have a computational complexity that depends on the implementation of the trees.

As a first approximation, however, we can consider the cost of these operations to be a unit. Then, the computational complexity of the algorithm is bounded by the number of such operations in a run. A brief inspection of Algorithm~\ref{algo:rec} reveals that the complexity of a run is proportional to the sum of the number of nodes of $D_t^n$ for $t$ from $\tau(n)$ to $-1$. Taking into account the complexity of the basic tree operations, this would typically lead to a complexity $O\left(\sum_{t=\tau(n)}^{-1} |D_t^n| \log |D_t^n|\right)$.

Thus, an analysis of the computational complexity  of Algorithm~\ref{algo:rec} requires simultaneous bounding of the number of iterations $\tau(n)$ and the size of the trees $D_t^n$. For a general kernel $P$, this involves not only the mixing properties of the corresponding process, but also the oscillation of the kernel itself, which amounts to a very challenging task that surpasses the scope of this paper. Nevertheless, Section~5 contains some analytical elements. The problematic issues are considered successively: First we give a crude bound on $\tau(n)$, then we prove a bound on the size of $D_t^n$s for finite memory processes.

\section{Bounding the size of $D_t^n$}\label{sec:boundingSize}
 We establish sufficient conditions for the algorithm to terminate and define the bound on the expectation of the depth of $D_t^n$. We then focus on the special, yet important, case of (finite) variable length Markov chains.

\subsection{Almost sure termination of the coupling scheme}
In general, the size of the CSD $D_t^n$ may be arbitrary with a positive probability.
The conditions that ensure the finiteness of $\tau(n)$, defined above, from which the bounds on $D_t^n$ can be deduced are given in \cite{CometsFernandezFerrari02longMem}.
However, these conditions are quite restrictive: in particular, they require that $A_0(\varepsilon)>0$.
The hybrid simulation scheme used in \cite{DesantsiPiccioni12perfectSim}, which allows for $A_0(\varepsilon)=0$, somewhat relaxes these conditions.

We define a crude bound, ignoring the coalescence possibilities of the algorithm: the depth of the current tree at time $t$ is defined as $L^n_t = \depth(D_t^n)$. Thereby an immediate inspection of Algorithm~\ref{algo:rec}  yields
\[\begin{cases}
L_t^t \leq \max\{X_t, L_{t+1}^{t+1}-1, 1\} & \hbox{ if } t<n \hbox{\;, and}\\
L_t^n \leq \max\{X_t, L_{t+1}^n-1\} & \hbox{ if }  t\geq n \;,\\

\end{cases}\]
where $X_t = \depth(D(U_t))$ represents i.i.d. random variables such that, for all $k\in \N,\, \P(X_t\leq k) = A^-_k$.
Therefore, \[\P(L_t^n\leq k) \geq \P(L_{t+1}^n \leq k+1) A^-_k \geq  \prod_{j=k}^{k-t-1} A^-_j \;. \]
As previously demonstrated in \cite{CometsFernandezFerrari02longMem}, it results that $\tau(n)$ is almost surely finite as soon as 
\[\sum_{m\geq 0}^\infty\prod_{k=0}^m A_k^- = \infty\;.\]
In addition, it follows that
\[\E[L_t^n] \leq \sum_{k=1}^n \left( 1-\prod_{j=k}^\infty A^-_j\right) \;.\]

\subsection{Finite context trees}\label{subsec:finiteCT}
It is easy to upper-bound the size of $D^n_t$ independently of $t\leq n$ for at least one case: when the kernel $P$ actually defines a \emph{finite context tree}, that is, when the mapping $\w \to P(\cdot |\w)$ is piecewise constant.
In other words, $P(\cdot | s)$ is a singleton for each $s\in D$, where $D$ denotes  the minimal CSD of this mapping.

Yet, the simulation scheme described above is useful even in that case: Although the ``plain'' Propp-Wilson algorithm could be applied to the first-order Markov chain $(X_{t+1:t+\depth(D)})_{t\in\Z}$ on the extended state space $G^{\depth(D)}$, the computational complexity of such an algorithm might become rapidly intractable if the depth $\depth(D)$ is large. In contrast, the following property shows that our algorithm retains a possibly much more limited complexity.
It is precisely because of these qualities of \emph{parsimony} that finite context trees have proved successful in many applications, from information theory and universal coding (see \cite{rissanen1983,willems1995,csiszar2006,garivier2006}) to biology (\cite{bejerano2001a,busch2009}) and linguistics~\cite{galves2009}.

We say that a CSD $D$ is \emph{prefix-closed} if every prefix of any sequence in $D$ is the suffix of an element of $D$:
\[\forall s\in D, \forall  k\leq |s|, \exists w\in D: w\hassuff s_{-|s|:-k}  \;.\]
A prefix-closed CSD satisfies the following property:
\begin{lemma}
If $D$ is a prefix-closed CSD, then $ha \hassuff D$ for all $h\in D$ (or, equivalently, for all $h\hassuff D$) and for all $a\in G$.
\end{lemma}
\textbf{Proof:}
If $h\in G^*$ is such that $ha \hassuff D$ does not hold for some $a\in G$, then (because D is a CSD) there exists $s\in D$ and $s' \in G^*\backslash \{\epsilon\}$ such that $s=s'ha$. But then $s'h$ is a prefix of $s$ and, by the prefix-closure property, there exists  $w \in D$ such that  $w\hassuff s'h$. Thus, we cannot have $h\hassuff D$.

We define the \emph{prefix closure} $\pcD$ of a CSD $D$ as the minimal prefix-closed set containing $D$, that is, the set of maximal elements (for the partial order $\hassuff$) of
\[\tilde{D} = \left\{ s_{-|s|:-k} : s\in D, \; k \leq |s|  \right\}\;.\]
In other words, $\pcD$ is the smallest set such that, for all $w\in \tilde{D}$, there exists $s\in \pcD$ such that $s\hassuff w$.

Clearly, $|\pcD|\leq |\tilde{D}| \leq |D| \times \depth(D)$. This bound is generally pessimistic: many CSDs are already prefix-closed, and for most CSDs, $|\pcD|$ is of the same order of magnitude as $|D|$. But, in fact, for each positive integer $n$, we can show that there exists a CSD $D$ of size $n$ such that $|\pcD| \geq c |D|^2$ for some constant $c\approx 0.4$.

Now, assume that $D\neq\{\epsilon\}$, i.e.
that $P$ is not memoryless.
\begin{prop}\label{prop:majD}
For each $t\leq n$ and for all $k<t, \pcD\hassuff D_t^n$.  Therefore,
\[|D_t^n|\leq \left|\pcD\right| \leq |D| \times \depth(D)\;.\]
\end{prop}

\textbf{Proof:} We show that $\pcD\hassuff D_t^t$ by induction on $t$. First, as $P$ is not memoryless, $\pcD\hassuff  D_{-1}^{-1} = G$.
Second, assume that $\pcD\hassuff  D_{t+1}^n$: it is sufficient to prove that $H_{t}^n$ (or $H_t^t$, if $t\geq n$) is $\pcD$-constant.
Observing that $\pcD\hassuff D$, for every $U_t\in[0,1[$ and for every $s\in \pcD$, it holds that $\phi_P(U_t, s)$ is a singleton $\{g_t[s]\}$.
By successively using the lemma and the induction hypothesis, we have $sg_t[s]\hassuff \pcD\hassuff D_t^n$, thus $H_t^n(s) = H_{t+1}^n(sg_t[s])$ is also a singleton.
Finally, for $t=n$, $H_t^t$ is $\pcD$-constant because $\pcD\neq\{\epsilon\}$.

It has to be emphasized that Proposition~\ref{prop:majD} provides only an upper-bound on the size of $D_t^n$: in practice, $D_t^n$ is often observed to be much smaller. Even for non-sparse, large order Markov chains of an order of $d$, it is possible for Algorithm~\ref{algo:rec} to be faster than the Propp-Wilson algorithm on $G^d$, which generally requires the consideration of $|G|^d$ states at each iteration. Interested readers may want to run the matlab experiments available at \url{http://www.math.univ-toulouse.fr/~agarivie/Telecom/context/}.

\section{Example: A continuous process with long memory}\label{sec:examples}
This section briefly illustrates the strengths of Algorithm~\ref{algo:rec} in comparison with the other existing CFTP algorithms for infinite memory processes. We focus on a process that cannot be simulated by other methods, although, of course, Algorithm~\ref{algo:rec} is also relevant for all processes mentioned in \cite{CometsFernandezFerrari02longMem, DesantsiPiccioni12perfectSim}, which we refer to for further examples.

\begin{figure}
\begin{center}
\includegraphics[width=10cm]{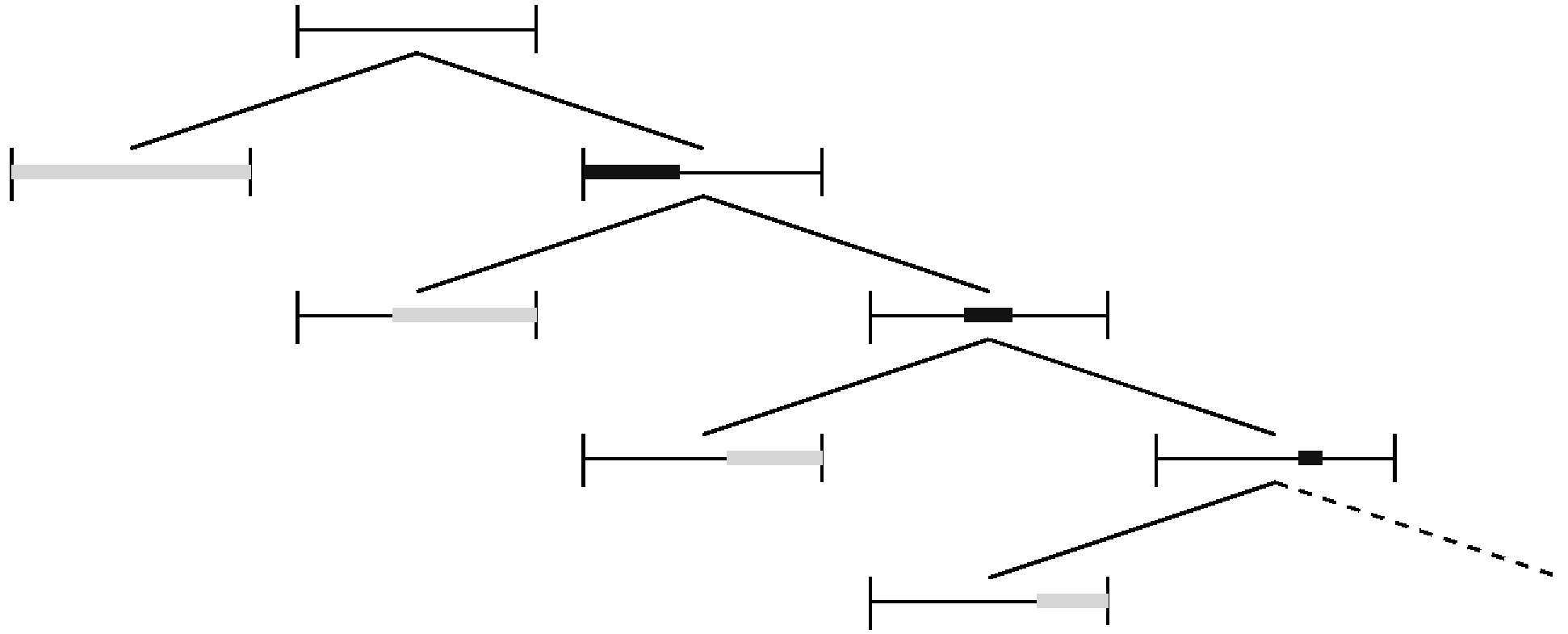}
\caption{Graphical representation of the update rule for the kernel defined in~\eqref{eq:defP}. Dark grey corresponds to $0$. Light grey corresponds to $1$.}
\label{fig:exCoupling}
\end{center}
\end{figure}

The example we consider involves a non-regenerating kernel on the binary alphabet $G=\{0,1\}$. It is such that $a_0=0$ and that the convergence of the coupling coefficients is slow, so that neither the perfect simulation scheme of
 \cite{CometsFernandezFerrari02longMem}, nor its improvement by \cite{DesantsiPiccioni12perfectSim} can be applied. Yet, a probabilistic upper-bound on the stopping time $\tau$ of Algorithm~\ref{algo:rec} can be given, which proves that there exists a compatible stationary process.
For all $k\geq 0$, let 
\begin{equation}\label{eq:defP}
 P(0|01^k) = 1-1/\sqrt{k} \;.
\end{equation}
Figure~\ref{fig:exCoupling} shows the coupling coefficients of $P$.
As $P(1|0) = \lim_{k\to \infty} P(0|01^k) = 1$,  $a_0=0$.
Moreover, for $k\geq 0$, it holds that $A_{k+1} = A_k(01^{k}) = 1-1/\sqrt{k}$, so that
\[\sum_n \prod_{k=2}^n A^-_k <\infty\] and the continuity conditions of \cite{CometsFernandezFerrari02longMem,DesantsiPiccioni12perfectSim} do not apply.

We demonstrate that the algorithm described above can be used to simulate samples of a process $X$ with the specification $P$ (so that, in particular, such a process exists; uniqueness is straightforward).
It is sufficient to show that the stopping time $\tau(1)$ is almost surely finite.
In fact, $-\tau(1)$ is stochastically upper-bounded three times by a geometric variable of parameter $3/(2\sqrt{2})-1$.
To simplify notations, we write $ H_{t} = H_t^{-1}$ and $\0=(\dots,0,0)\in\G$.
For every $t<-2$, if $U_{t-1}\leq 1-1/\sqrt{2}$, if $U_{t}>1-1/\sqrt{2}$,  and if $U_{t+1} \leq 1-1/\sqrt{2}$, then, for every $\w\in\G$, we can see that
\begin{itemize}
 \item $U_{t+1} \leq 1-1/\sqrt{2}$ implies that $f_{t+1}(\w 1)=\w 10$ and $H_{t+1}(\w 1) = H_{t+2}(\0)$;
 \item $U_{t}>1-1/\sqrt{2}$  implies that $f_{t}(\w01)=\w 011$ and $H_{t}(\w01)=H_{t+1}(\w011)=H_{t+2}(\0)$, whereas $f_{t}(\w0)=\w 01$ and  $H_{t}(\w0)=H_{t+1}(\w01)=H_{t+2}(\0)$;
 \item $U_{t-1} \leq 1-1/\sqrt{2}$ implies that $f_{t-1}(\w 1) = \w 10$ and $f_{t-1}(\w 0) = \w01$, so that $H_{t-1}(\w 0) = H_{t-1}(\w 1) = H_{t+2}(\0)$ and $\tau \geq t-1$.
\end{itemize}
\begin{figure}
\begin{center}
\input{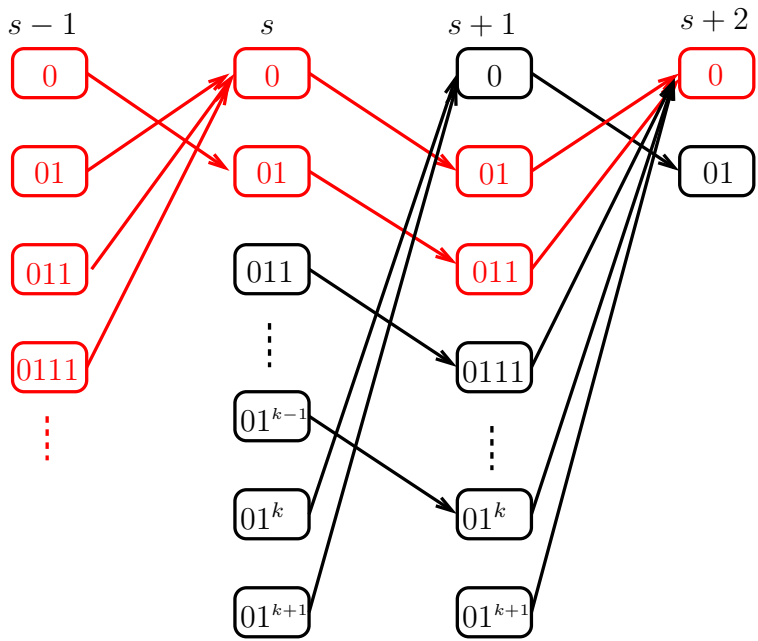}
\caption{Convergence of the simulation scheme.}
\label{fig:exampleCV}
\end{center}
\end{figure}
For every negative integer $k$, let $E_k = \{U_{3k-1}\leq 1-1/\sqrt{2}\}\cap\{U_{3k}>1-1/\sqrt{2} \}\cap\{ U_{3k+1} \leq 1-1/\sqrt{2}\}$. The events $(E_k)_{k<0}$ are independent and of probability $3/(2\sqrt{2})-1$, which gives the result.

Thus, the algorithm converges fast. However, the dictionaries involved in the simulation can be very large. In fact, it is easy to see  that the depth $X_t$ of $D(U_t)$ has no expectation: $\P(X_t\geq k) = 1/\sqrt{k}$.
Of course, because of the very special shape of $D_t^n$, ad hoc modifications of the algorithm would allow us to easily draw arbitrary long samples with low computational complexity. Moreover, the paths of the renewal processes can be simulated directly. Nevertheless, this example illustrates the weakness of the conditions required by Algorithm~\ref{algo:rec}; it also shows that neither regeneration nor a rapid decrease of the coupling coefficients are necessary conditions for perfect simulation. It is easy to imagine more complicated variants for which no other sampling method is currently known.

To conclude, a simple modification of this example shows that continuity is absolutely not necessary to ensure convergence because the proof also applies to any kernel $P'$, such that $P'(0|01^k) = 1-1/\sqrt{k}$ for $0\leq k\leq 2$, and $P'(0|\w11) \geq 1-1/\sqrt{2}$ for any $\w\in \G$.
Gallo, who has studied such a phenomenon in \cite{Gallo10perfectSim}, gives sufficient conditions on the shape of the trees (together with bounds on transition probabilities) to ensure convergence of his coupling scheme. His approach is quite different and does not cover the examples presented here.

\section*{Acknowledgments}
We would like to thank the reviewers for their help in improving the redaction of this paper and for pointing out Kendall's `coupling from and into the past' method (named by Wilson).
We further thank Sandro Gallo, Antonio Galves, and Florencia Leonardi (Numec, Sao Paulo) for the stimulating discussions on chains of infinite memory. This work received support from USP-COFECUB (grant 2009.1.820.45.8).

\bibliographystyle{plain}
\nocite{huber06fastPS,kendall98PS,wilson00PSreadonce}
\bibliography{biblio}

\section*{Appendix}

\subsection*{Proof of Lemma~\ref{lem:phiPC}}
Let $u\in[0,1[$. The uniform continuity of $P$ implies that there exists $\epsilon$ such that, if $\delta(\w,\z)\leq \epsilon$, then $\big|P(\cdot|\w) - P(\cdot|\z) \big|_{TV}<(1-u)/|G|$.
However, Equation~\eqref{eq:coupling} implies that $\phi_P(u, \w) = \phi_P(u, \z)$.

\subsection*{Proof of Proposition~\ref{prop:etavsA}}
For the upper-bound, observe that
\begin{align*}
\eta_P(s) &= \sup\left\{\big|P(\cdot|\w) - P(\cdot|\z) \big|_{TV}  : \w, \z \in \tree(s) \right\}\\
 &= \sup\left\{1-\sum_{a\in G}P(a|\w) \wedge P(a|\z) : \w, \z \in \tree(s) \right\}\\
 &= 1 - \inf\left\{\sum_{a\in G}P(a|\w) \wedge P(a|\z) : \w, \z \in \tree(s) \right\}\\
 &\leq 1 - \sum_{a\in G}\inf\left\{P(a|\w)\wedge P(a|\z) : \w,\z \in \tree(s) \right\}\\
 &= 1 - \sum_{a\in G}\inf\left\{P(a|\w) : \w \in \tree(s) \right\}\\
 &= 1 - A_{|s|}(s)\;.
\end{align*}
For the lower-bound, let $\epsilon>0$, let $\w\in\tree(s)$ and $b\in G$ be such that, for all $\z\in \tree(s)$ and  for all $a\in G, P(a|\z)\geq P(b|\w)-\epsilon$. Then, for all $z\in\G$ and all $a\neq b, P(a|\z)\geq P(a|\w) - \eta_P(s)$. Thus,
\begin{align*}
&A_{|s|}(s) = \sum_{a\in G} P(a|\w) + \inf\left\{P(a|\z) - P(a|\w): \z\in \tree(s) \right\} \\
&\geq  1 + \inf\left\{P(b|\z) - P(b|\w): \z\in \tree(s) \right\}  + \sum_{a \neq b} \inf\left\{P(a|\z) - P(a|\w): \z\in \tree(s) \right\} \\
&= 1 - \epsilon - (|G|-1)\eta_P(s) \;.
\end{align*}
Because $\epsilon$ is arbitrary, the result follows.

\subsection*{Proof of Proposition~\ref{prop:Pcont}}
The equivalence of (i) and (ii) is obvious by definition. 
The equivalence with (iii) is a simple consequence of Proposition~\ref{prop:etavsA}.
Similarly, (iii) follows from (i): if $P$ is continuous on the compact set  $\G$, then it is uniformly continuous, and
\[\varphi(k) = \sup_{s\in G^{-k}} \eta_P(s) \to 0\]
as $k$ goes to infinity. But, by Proposition \ref{prop:etavsA}, $A^-_k \geq 1-|G| \varphi(k)$.
Finally, (iii) implies (ii).

The equivalence of (ii) and (iii) can also be proved as a consequence of Dini's theorem (see\cite{rudin76principles}, Theorem 7.13, page 150): based on the definition $\tilde{A}_k(\w) = A_k(w_{-k:-1})$, the sequence $\left(\tilde{A}_k\right)_k$ is an increasing sequence of continuous functions that converges pointwise to the (continuous) constant function $1$, thus the convergence is uniform.

\subsection*{Proof of Proposition~\ref{prop:partition}}
Let $m$ (resp. $M$) be the minimal (resp. maximal) element of $G$.
Then, for any integer $k$, $\alpha_k(m | w_{-k:-1}) = A_{k-1}(w_{-k+1:-1})$, $\beta_k(M | w_{-k:-1}) = A_k(w_{-k:-1})$, and
\[\left[A_{k-1}(w_{-k+1:-1}), A_k(w_{-k:-1})\right[ = \bigsqcup_{g\in G}\left[\alpha_k(g | w_{-k:-1}), \beta_k(g | w_{-k:-1})\right[ \;.\]
But $A_{-1}(\varepsilon)=0$, and thus the result follows from the continuity assumption: $A_k(w_{-k:-1}) \to 1$ when $k$ goes to infinity.

\subsection*{Proof of Proposition~\ref{prop:phiOk}}
We need to prove that if $U\sim\U([0,1[)$, then for all $\w\in \G$ the random variable $\phi_P(U, \w)$ has a distribution $P(\cdot | \w)$.
It is sufficient to prove that, for all $g\in G$, \[\sum_{l=0}^\infty \beta_l(g | w_{-l:-1}) - \alpha_l(g | w_{-l:-1}) = P(g| \w).\]
For any integer $k$, it holds that 
\begin{align*}
 \sum_{l=0}^k \beta_l(g | w_{-l:-1}) - \alpha_l(g | w_{-l:-1}) &= \sum_{l=0}^k a_{l}(g | w_{-l:-1}) - a_{l-1}(g | w_{-l+1:-1}) \\
 &=  a_{k}(g | w_{-k:-1})\;.
\end{align*}
As an increasing sequence upper-bounded by $P(g|\w)$, $a_k(g|w_{-k:-1})$ has a limit $Q(g|\w)\leq P(g|\w)$ when $k$ tends to infinity.
By continuity,
\begin{multline*}\sum_{g\in G }  Q(g|\w)  = \sum_{g\in G } \lim_{k\to\infty} a_k(g|w_{-k:-1}) \\ = \lim_{k\to\infty}  \sum_{g\in G}a_k(g|w_{-k:-1}) = \lim_{k\to\infty} A_k(w_{-k:-1})=  1,\end{multline*}
and, because $\sum_{g\in G}P(g|\w)=1$, this implies that, for all $g\in G$, $Q(g|\w) = P(g|\w)$.

\noindent The last part of the proposition is immediate: for $\w$ and $\z\in \tree(s)$ and for all $k\leq |s|, w_{-k:-1} = z_{-k:-1}$, and
\begin{multline*}
 \bigsqcup_{g\in G, k\leq|s|}\left[\alpha_k(g|w_{-k:-1}), \beta_k(g|w_{-k:-1})\right[ \\=\bigsqcup_{g\in G, k\leq|s|}\left[\alpha_k(g|z_{-k:-1}), \beta_k(g|z_{-k:-1})\right[ = [0, A_{|s|}(s)[ \;.
\end{multline*}

\end{document}